\documentclass[twoside,12pt]{amsart}
\usepackage{mathrsfs}
\usepackage{amsfonts}
\newtheorem{thm}{Theorem}[section]
\newtheorem{cor}[thm]{Corollary}
\newtheorem{lem}[thm]{Lemma}
\newtheorem{prop}[thm]{Proposition}

\makeatletter
\def\serieslogo@{}
\makeatother
\makeatletter
\def\@setcopyright{}
\makeatother

\theoremstyle{definition}

\theoremstyle{remark}
\newtheorem{rem}{Remark}[section]

\begin{document}

\title{Weak entropy solutions of nonlinear reaction-hyperbolic
systems for axonal transport}

\author{Hao YAN}
\address{Zhou Pei-Yuan Center for Appl. Math.\\
Tsinghua University\\
Beijing 100084, China} \email{yanhao06@mails.tsinghua.edu.cn}

\author{Wen-An YONG}
\address{Zhou Pei-Yuan Center for Appl. Math.\\
Tsinghua University\\
Beijing 100084, China} \email{wayong@tsinghua.edu.cn}




\begin{abstract}
This paper is concerned with a class of nonlinear
reaction-hyperbolic systems as models for axonal transport in
neuroscience. We show the global existence of entropy-satisfying
BV-solutions to the initial-value problems by using hyperbolic-type
methods. Moreover, we rigorously justify the limit as the
biochemical processes are much faster than the transport ones.

\vskip 10 pt

\end{abstract}

\maketitle \markboth{Hao YAN and W.-A. YONG}{Weak entropy solutions
of nonlinear reaction-hyperbolic systems for axonal transport}

\tableofcontents

\section{Introduction}

The axonal transport is important for the maintenance and functions
of nerve cells. These cells are also called {\it neurons}.
A neuron consists of three parts mainly: cell body, dendrites and a
single axon. The axon is a long and thin pipe whose length can
exceed 10,000 times its diameter. It is this axon that distinguishes
neurons from other cells. The axon is responsible for signal
transmission in the nervous system. Its cytoplasm does not contain
rough endoplasmic reticulum and therefore its proteins can only be
transported from the cell body, where all proteins are synthesized.

The transport proceeds as follows. Proteins are stored in vesicles
as cargos. The vesicles are attached to kinesin (anterograde motors)
or dynein (retrograde motors) proteins. These motor proteins drive
the vesicles to walk along the cytoskeletal microtubules as track.
Here the kinesin proteins move the vesicles from the cell body to
synapse (anterograde transport), while the dynein proteins move the
vesicles in the opposite direction (retrograde transport). During
the transport, many biochemical processes are possible. For example,
the cargos can leave its track, can switch its motor proteins from
kinesin to dynein or vice verse, and can move back onto the track.
Thus, we can divide the cargos into a number of subpopulations, such
as free vesicles,
vesicle-kinesin compounds off track, moving vesicle-dynein compounds
on track, etc.

As the axon is long and thin, it is reasonable to assume the
transport only along the longitudinal direction of the axon. Denote
by $x>0$ the distance down the axon from the cell body
which is located at $x = 0$. Let $u_i = u_i(x, t)$ be the
concentration at space-time $(x, t)$ of the $i$-th subpopulations.
According to Reed and Blum \cite{RB}, the mathematical model for
axonal transport is partial differential equations of the form:
\begin{equation*}\label{pde0}
\partial_t u_i + \lambda_i\partial_x u_i = F_i (u_1, u_2, \cdots, u_r), \qquad i =
1, 2, \cdots, r .
\end{equation*}
Here the term $\lambda_i\partial_x u_i$ accounts for the transport
of the $i$-th subpopulation with constant velocity $\lambda_i$, and
$F_i (u_1, u_2, \cdots, u_r)$ describes the biochemical processes of
the constituents. It
is well recognized that the biochemical processes are much faster
than the transport in biosystems.


In \cite{CBF, FC, FH, RVB}, the authors studied the linear case,
where $F_i(u_1, u_2, \cdots, u_r)$ is linear with respect to the
$u_j$'s, in order to explain the approximate traveling waves
observed in experiments. Especially, in \cite{FH} Friedman and Hu
used parabolic-type estimates to analyse the diffusive limit of the
linear systems. However, it seems uneasy to deal with the nonlinear
problems with the parabolic-type techniques. On the other hand, in
\cite{Ca} Carr showed the existence of global classical solutions to
a class of nonlinear models with source terms of the form
$$
F_i(u_1, u_2, \cdots, u_r) = \begin{cases} f_1 , & i = 1\\
f_i - f_{i - 1} , & 1 < i < r \\
-f_{r-1} , & i = r
\end{cases},
$$
where $f_i = f_i(u_i, u_{i+1})$ is a continuously differentiable
function of two variables. We notice that in applications $f_i$ is a
polynomial of $u_i$ and $u_{i+1}$. In addition, the model in
\cite{CBF} is an example where $f_i$ depends only on $u_i$ and
$u_{i+1}$.

In this paper, we consider the same nonlinear systems, as in
\cite{Ca}, but with a small parameter $\epsilon > 0$:
 \begin{eqnarray}\label{pde1}
\left\{ \begin{array}{ll}
(\partial_t+\lambda_{1}\partial_x)u_1(x,t) & = \frac{1}{\epsilon}f_{1}(u_1,u_2), \\
\qquad\cdots  \\
(\partial_t+\lambda_{i}\partial_x)u_i(x,t) & =
-\frac{1}{\epsilon}f_{i-1}(u_{i-1},u_{i})+\frac{1}{\epsilon}f_{i}(u_i,u_{i+1}), \\
\qquad \cdots  \\
(\partial_t+\lambda_{r}\partial_x)u_r(x,t) & =
-\frac{1}{\epsilon}f_{r-1}(u_{r-1},u_{r}) .
\end{array} \right.
\end{eqnarray}
Here the small parameter $\epsilon$ characterizes the fact that the
biochemical processes are much faster than the transport. We assume,
throughout this paper, that each $f_i$ is strictly decreasing with
respect to the first argument and strictly increasing with respect
to the second. In addition, we assume that $f_i(0,0)=0$ and, for
fixed $v$, there exists $w$ such that $vf_i(v, w) \geq 0$. These
assumptions are consistent with those used in \cite{Ca, CBF, FC, FH,
RB, RVB}.

We will regard (\ref{pde1}) as a hyperbolic relaxation system
\cite{Y2} and use the corresponding techniques to study it. In
particular, we will use a difference scheme to show the global
existence of entropy-satisfying BV-solutions $(u_1^\epsilon,
u_2^\epsilon, \cdots, u_r^\epsilon)$ to the initial-value problems
of (\ref{pde1}) and investigate the limit as $\epsilon$ goes to
zero. For $r\leq 3$ and for linear problems, we prove that the limit
$(u_1^0, u_2^0, \cdots, u_r^0)$ is also an entropy-satisfying
BV-solution to the so-called equilibrium or reduced system of
$(\ref{pde1})$:
\begin{eqnarray}\label{pde2}
\left\{\begin{array}{ll}
\partial_t(u_1+u_2+\cdots+u_r) +
\partial_x(\lambda_1u_1+\lambda_2u_2+\cdots+\lambda_ru_r)=0,\\
f_1(u_1,u_2) = f_2(u_2,u_3) = \cdots = f_{r-1}(u_{r-1},u_r)=0 .
\end{array} \right.
\end{eqnarray}

Although the $BV$-framework is quite standard for nonlinear
hyperbolic problems (see, e.g.,  \cite{Da, RY, Y1}), some innovatory
ideas seem necessary to carry our the details. In particular, we use
the Brouwer fixed-point theorem, the $BV$-estimate and the special
structure of the source terms in (\ref{pde1}) to derive the
existence, boundedness and time-Lipschitz continuity of the
difference solutions. It seems not so easy to obtain the
time-Lipschitz continuity! Moreover, for an arbitrarily given convex
entropy function for the reduced system ({\ref{pde2}), we construct
a dissipative entropy function for the original system (\ref{pde1}).
In equilibrium, the constructed entropy function reduces to the
given one.

The paper is organized as follows. Section 2 is devoted to the
analysis of a difference scheme for the system $(\ref{pde1})$. The
time-Lipschitz continuity of the difference solutions is derived in
Section 3. In Section 4, we discuss entropy functions for the two
systems $(\ref{pde1})$ and $(\ref{pde2})$. The main results are
shown in Section 5.

\section{Difference Solutions}
\setcounter{equation}{0}

We begin with construction of approximation solutions to the
reaction-hyperbolic system $(\ref{pde1})$ by using difference
methods. For simplicity, we set
$$
\begin{array}{rl}
U = & (u_{1}, u_{2}, \cdots, u_{r})^{T}, \quad
\Lambda = \mbox{diag}(\lambda_1,\lambda_2,\cdots,\lambda_ r),
\\[4mm]
Q(U) = & (f_1, \cdots, -f_{i-1}+f_{i}, \cdots, - f_{r-1})^{T},
\end{array}
$$
where the superscript $T$ denotes the transpose of vectors or
matrices. Then the reaction-hyperbolic system $(\ref{pde1})$ can be
rewritten as
\begin{eqnarray}\label{pde3}
U_t + \Lambda U_x = \frac{1}{\epsilon}Q(U).
\end{eqnarray}
About this system, we make the following assumptions mentioned in
the introduction:\\

\begin{itemize}

\item[(1).] each $f_i$ is continuously differentiable, strictly decreasing
with respect to the first argument, and strictly increasing with
respect to the second;

\item[(2).] $f_i(0, 0) = 0$.\\

\end{itemize}
These assumptions are consistent with those used in \cite{Ca, CBF,
FC, FH, RB, RVB}.

Our difference approximation to $(\ref{pde3})$ is the following
semi-implicit upwind scheme
\begin{eqnarray}\label{dc1}
\frac{U_j^{n+1}-U_j^n}{\Delta
t}+\Lambda^+\frac{U_j^{n}-U_{j-1}^{n}}{\Delta
x}+\Lambda^-\frac{U_{j+1}^{n}-U_j^{n}}{\Delta x}=\frac{1}{\epsilon}
Q_j^{n+1}.
\end{eqnarray}
Here $\Delta x$ and $\Delta t$ denote the increments respectively in
$x$ and $t$; $U_{j}^n$ denotes the approximation of $U(x,t)$ over
the grid block $[x_j,x_{j+1})\times [t_n,t_{n+1})$ with $x_j=j\Delta
x$ and $t_n=n\Delta t$; $j=0,\pm1,\pm2,\cdots$ and $n=0,1,2,\cdots$;
$Q_j^n=Q(U_j^n);
\Lambda^+=diag(\lambda_1^+,\lambda_2^+,\cdots,\lambda_
r^+)\quad\mbox{and}\quad
\Lambda^-=diag(\lambda_1^-,\lambda_2^-,\cdots,\lambda_ r^-)$ with
\begin{eqnarray}\label{UPwind}
\lambda_{i}^+=\frac{\lambda_{i}+|\lambda_i|}{2} \ge 0
\quad\mbox{and}\quad \lambda_{i}^-=\frac{\lambda_{i}-|\lambda_i|}{2}
\le 0.
\end{eqnarray}
For $n=0$, we take
\begin{eqnarray}\label{inital}
U_j^0 = \begin{cases}
\frac{1}{\Delta x}\int_{x_j}^{x_{j+1}} U_0(x) dx, &
\mbox{if} \quad j \leq \frac{1}{\Delta x} \\[4mm]
0, & \mbox{if} \quad j > \frac{1}{\Delta x} \end{cases},
\end{eqnarray}
where $U_{0}(x):=(u_{10}(x),u_{20}(x),...,u_{r0}(x))^T$ is a bounded
measurable function of $x\in\mathbf{R}$. In addition, the 0 in
(\ref{inital}) can be replaced with any constant vector.

Throughout this paper, we assume that the grid sizes satisfy the
CFL-condition
\begin{eqnarray}\label{CFL}
\frac{\Delta t}{\Delta x}\max_i |\lambda_{i}|\le1.
\end{eqnarray}

To analyse the above scheme, we start with the following elementary
fact.

\begin{lem}\label{lem1}
Let $A=(a_{kl})_{n\times n}$ be a matrix satisfying $\sum_{k=1}^n
a_{kl}=1$ for $l=1,2,\cdots,n$ and $a_{kl}\leq 0 \ \textrm{for }
k\ne l$. Then $A$ is invertible, the spectral radius of $A^{-1}$ is
not bigger than 1, and the $1$-norm of $A^{-1}$ is 1.
\end{lem}

\begin{proof}
Let $\lambda$ be an eigenvalue of $A$. By the Gershgorin circle
theorem there is an integer $l$ such that
\begin{equation*}
|\lambda - a_{ll}|\leq\sum_{k\neq l}|a_{kl}| = a_{ll} - 1,
\end{equation*}
for $a_{ll}=1-\sum_{k\ne l} a_{kl}$ and $a_{kl}\leq 0 (k\ne l)$.
From this we easily deduce that $|\lambda|\geq1$. Thus $A$ is
invertible and the spectral radius of $A^{-1}$ is not bigger than 1.

Furthermore, it follows from $\sum_{k=1}^n a_{kl} = 1$ that
\begin{equation*}
(1,1,\cdots,1)A=(1,1,\cdots,1)
\end{equation*}
and thereby
\begin{equation*}\label{matrix1}
(1,1,\cdots,1)A^{-1}=(1,1,\cdots,1).
\end{equation*}
This shows that the sum of each column of $A^{-1}$ is also 1. On the
other hand, we set $D = \mbox{diag}(a_{11}, a_{22}, \cdots, a_{nn})$
and $F = D - A$. It is obvious that the elements of $FD^{-1}$ are
nonnegative,
\begin{equation*}
A=D-F=(I-FD^{-1})D
\end{equation*}
and thereby
\begin{equation*}
A^{-1}=D^{-1}(I-FD^{-1})^{-1} .
\end{equation*}
Moreover, with the Gershgorin circle theorem it is not difficult to
see that the spectral radius of $FD^{-1}$ is less than 1. Thus we
can write
\begin{equation*}
A^{-1}=D^{-1}(I+FD^{-1}+(FD^{-1})^2+\cdots),
\end{equation*}
which shows that the elements of $A^{-1}$ are nonnegative. Hence the
$1$-norm of $A^{-1}$ is 1 for the sum of each column of $A^{-1}$ is
1. This completes the proof.
\end{proof}

Set $V=(v_1, v_2, \cdots, v_r)^T, W=(w_1, w_2, \cdots, w_r)^T, $
\begin{eqnarray}
A_i& = &A_i(V,W) = \int_0^1 \frac{\partial f_i}{\partial
u_i}(\theta v_i + (1-\theta )w_i, \theta v_{i+1} + (1-\theta)w_{i+1})d\theta,\label{int1}\\
B_i& = &B_i(V,W) = \int_0^1 \frac{\partial f_i}{\partial
u_{i+1}}(\theta v_i + (1-\theta)w_i, \theta v_{i+1} +
(1-\theta)w_{i+1})d\theta.\label{int2}
\end{eqnarray}
Assumption (1) implies that $A_i<0$ and $B_i>0.$ By the mean-value
theorem we have
\begin{equation}
Q(V)-Q(W)=\mathcal {Q}(V,W)(V-W).
\end{equation}
From the definition of $Q$, we have
\begin{eqnarray*}\label{Jocobi}
\mathcal {Q}(V,W)=\left( \begin{array}{cccccc}
A_{1} & B_{1}  &   &  &    \\
-A_{1} & -B_{1}+A_{2}& B_{2}  & & \\
{}&  -A_{2}& -B_{2}+A_{3}& B_{3}& \\
{}& &  \ldots& \ldots &\ldots\\
 {}& &  & -A_{r-1}& -B_{r-1}
\end{array} \right).
\end{eqnarray*}
It is easy to verify that $I-\frac{\Delta t}{\epsilon}\mathcal
{Q}(V,W)$ satisfies the conditions in Lemma $\ref{lem1}.$ Therefore
we have

\begin{cor}\label{cor}
$I-\frac{\Delta t}{\epsilon}\mathcal {Q}(V,W)$ is invertible, the
spectral radius of $(I-\frac{\Delta t}{\epsilon}\mathcal
{Q}(V,W))^{-1}$ is not bigger than 1, and the $1$-norm of
$(I-\frac{\Delta t}{\epsilon}\mathcal {Q}(V,W))^{-1}$ is 1 for any
$V,W \in \mathbb{R}^r$.
\end{cor}

With these preparations, we study the difference scheme
$(\ref{dc1})$. In what follows, we will denote by $|\cdot|$ the
$1$-norm of vectors and matrices.
\begin{lem}\label{lem2}
Given $U_j^n$, $U_j^{n+1}$ can be uniquely determined by solving the
nonlinear algebraic equation $(\ref{dc1})$.
\end{lem}

\begin{proof}
From $f_i(0,0)=0$ we see that $Q(0)=0$ and
\begin{eqnarray*}
Q_j^{n+1}=Q_j^{n+1}-Q(0)=\mathcal {Q}(U_j^{n+1},0)U_j^{n+1}.
\end{eqnarray*}
Thus, we deduce from the difference scheme $(\ref{dc1})$ that
\begin{eqnarray*}
(I-\frac{\Delta t}{\epsilon}\mathcal
{Q}(U_j^{n+1},0))U_j^{n+1}=\frac{\Delta t}{\Delta
x}\Lambda^+U_{j-1}^n+(I-\frac{\Delta t}{\Delta
x}\Lambda^++\frac{\Delta t}{\Delta
x}\Lambda^-)U_{j}^n+(-\frac{\Delta t}{\Delta x}\Lambda^-)U_{j+1}^n.
\end{eqnarray*}

With this equation, we construct a mapping $F:V\rightarrow W$ as
follows
\begin{eqnarray*}
(I-\frac{\Delta t}{\epsilon}\mathcal {Q}(V,0))W=\frac{\Delta
t}{\Delta x}\Lambda^+U_{j-1}^n+(I-\frac{\Delta t}{\Delta
x}\Lambda^++\frac{\Delta t}{\Delta
x}\Lambda^-)U_{j}^n+(-\frac{\Delta t}{\Delta x}\Lambda^-)U_{j+1}^n.
\end{eqnarray*}
Due to Corollary $\ref{cor}$, we can write
\begin{eqnarray*}
W=(I-\frac{\Delta t}{\epsilon}\mathcal {Q}(V,0))^{-1}[\frac{\Delta
t}{\Delta x}\Lambda^+U_{j-1}^n+(I-\frac{\Delta t}{\Delta
x}\Lambda^++\frac{\Delta t}{\Delta
x}\Lambda^-)U_{j}^n+(-\frac{\Delta t}{\Delta x}\Lambda^-)U_{j+1}^n]
\end{eqnarray*}
and thereby
\begin{eqnarray*}
|W|&\leq& |(I-\frac{\Delta t}{\epsilon}\mathcal
{Q}(V,0))^{-1}||\frac{\Delta t}{\Delta
x}\Lambda^+U_{j-1}^n+(I-\frac{\Delta t}{\Delta
x}\Lambda^++\frac{\Delta t}{\Delta
x}\Lambda^-)U_{j}^n+(-\frac{\Delta t}{\Delta
x}\Lambda^-)U_{j+1}^n|\\&=& |\frac{\Delta t}{\Delta
x}\Lambda^+U_{j-1}^n+(I-\frac{\Delta t}{\Delta
x}\Lambda^++\frac{\Delta t}{\Delta
x}\Lambda^-)U_{j}^n+(-\frac{\Delta t}{\Delta x}\Lambda^-)U_{j+1}^n|.
\end{eqnarray*}
Therefore, $W$ is bounded for fixed $j$ and $n$. By the Brouwer
fixed-point theorem, the mapping $F$ has fixed points for $F$ is
continuous. This shows the existence of $U^{n+1}_j$.

To see the uniqueness, we assume that both $U_1$ and $U_2$ satisfy
the difference scheme $(\ref{dc1})$. With the mean-value theorem, we
see that
\begin{eqnarray*}
U_1-U_2=\frac{\Delta t}{\epsilon}(Q(U_1)-Q(U_2))=\frac{\Delta
t}{\epsilon}\mathcal {Q}(U_1,U_2)(U_1-U_2)
\end{eqnarray*}
and thereby
\begin{eqnarray*}
(I-\frac{\Delta t}{\epsilon}\mathcal {Q}(U_1,U_2))(U_1-U_2)=0.
\end{eqnarray*}
It follows from Corollary $\ref{cor}$ that $U_1=U_2$. This shows the
uniqueness and hence the proof is complete.
\end{proof}

\begin{rem}
We could use the contraction mapping principle to show the existence
and uniqueness of $U_j^{n+1}$. However, it requires that $\Delta t$
is much smaller than $\epsilon$. On the other hand, this lemma does
not tell the uniform boundedness of the difference solution $U_j^n$.
\end{rem}

Next we establish the $L^1$-stability of the difference scheme
$(\ref{dc1})$.
\begin{lem}\label{lem3}
Let $U_j^n$ and $V_j^n$ be two solutions of the difference scheme
$(\ref{dc1})$ with initial data $U_j^0$ and $V_j^0$, respectively.
Then it holds that
\begin{eqnarray*}
\sum_{j=-\infty}^{+\infty}|U_{j}^{n}-V_j^{n}|\leq\sum_{j=-\infty}^{+\infty}|U_{j}^{0}-V_j^{0}|
\end{eqnarray*}
for all $n\ge 0$.
\end{lem}

\begin{proof}
It follows from the difference scheme $(\ref{dc1})$ that
\begin{eqnarray*}
U_j^{n+1}-V_j^{n+1}&=&\frac{\Delta
t}{\epsilon}(Q(U_j^{n+1})-Q(V_j^{n+1}))+\frac{\Delta t}{\Delta
x}\Lambda^+(U_{j-1}^n-V_{j-1}^n)\\&{}&+(I-\frac{\Delta t}{\Delta
x}\Lambda^++\frac{\Delta t}{\Delta
x}\Lambda^-)(U_{j}^n-V_{j}^n)+(-\frac{\Delta t}{\Delta
x}\Lambda^-)(U_{j+1}^n-V_{j+1}^n).
\end{eqnarray*}
Thus we deduce from the mean-value theorem and Corollary $\ref{cor}$
that
\begin{eqnarray*}
(I-\frac{\Delta t}{\epsilon}\mathcal
{Q}(U_j^{n+1},V_j^{n+1}))(U_j^{n+1}-V_j^{n+1})=\frac{\Delta
t}{\Delta x}\Lambda^+(U_{j-1}^n-V_{j-1}^n)\\
+(I-\frac{\Delta t}{\Delta x}\Lambda^++\frac{\Delta t}{\Delta
x}\Lambda^-)(U_{j}^n-V_{j}^n)+(-\frac{\Delta t}{\Delta
x}\Lambda^-)(U_{j+1}^n-V_{j+1}^n) ,
\end{eqnarray*}
\begin{eqnarray*}
U_j^{n+1}-V_j^{n+1}&=&(I-\frac{\Delta t}{\epsilon}\mathcal
{Q}(U_j^{n+1},V_j^{n+1}))^{-1}[\frac{\Delta t}{\Delta
x}\Lambda^+(U_{j-1}^n-V_{j-1}^n)\\&{}&+(I-\frac{\Delta t}{\Delta
x}\Lambda^++\frac{\Delta t}{\Delta
x}\Lambda^-)(U_{j}^n-V_{j}^n)+(-\frac{\Delta t}{\Delta
x}\Lambda^-)(U_{j+1}^n-V_{j+1}^n)],
\end{eqnarray*}
\begin{eqnarray*}
|U_j^{n+1}-V_j^{n+1}|&\leq&|(I-\frac{\Delta t}{\epsilon}\mathcal
{Q}(U_j^{n+1},V_j^{n+1}))^{-1}||\frac{\Delta t}{\Delta
x}\Lambda^+(U_{j-1}^n-V_{j-1}^n)\\&{}&+(I-\frac{\Delta t}{\Delta
x}\Lambda^++\frac{\Delta t}{\Delta
x}\Lambda^-)(U_{j}^n-V_{j}^n)+(-\frac{\Delta t}{\Delta
x}\Lambda^-)(U_{j+1}^n-V_{j+1}^n)|\\
&\leq&|\frac{\Delta t}{\Delta
x}\Lambda^+(U_{j-1}^n-V_{j-1}^n)|\\&{}&+|(I-\frac{\Delta t}{\Delta
x}\Lambda^++\frac{\Delta t}{\Delta
x}\Lambda^-)(U_{j}^n-V_{j}^n)|+|(-\frac{\Delta t}{\Delta
x}\Lambda^-)(U_{j+1}^n-V_{j+1}^n)|,
\end{eqnarray*}
and thereby
\begin{eqnarray*}
&& \sum_{j=-\infty}^{+\infty}|U_j^{n+1}-V_j^{n+1}|
\leq\sum_{j=-\infty}^{+\infty}|\frac{\Delta t}{\Delta x}
\Lambda^+(U_{j-1}^n-V_{j-1}^n)|\\
&& \qquad + \sum_{j=-\infty}^{+\infty}|(I-\frac{\Delta t}{\Delta
x}\Lambda^++\frac{\Delta t}{\Delta
x}\Lambda^-)(U_{j}^n-V_{j}^n)|+\sum_{j=-\infty}^{+\infty}|(-\frac{\Delta
t}{\Delta x}\Lambda^-)(U_{j+1}^n-V_{j+1}^n)|\\
& = &\sum_{j=-\infty}^{+\infty}[|\frac{\Delta t}{\Delta
x}\Lambda^+(U_{j}^n-V_{j}^n)|+|(I-\frac{\Delta t}{\Delta
x}\Lambda^++\frac{\Delta t}{\Delta
x}\Lambda^-)(U_{j}^n-V_{j}^n)|\\
&& \qquad + |(-\frac{\Delta t}{\Delta
x}\Lambda^-)(U_{j}^n-V_{j}^n)|].
\end{eqnarray*}
In the last step we use the induction assumption that
$\sum_{j=-\infty}^{+\infty}|U_{j}^n-V_j^n|<\infty$, which is true if
$\sum_{j=-\infty}^{+\infty}|U_{j}^{0}-V_j^{0}|<\infty$. On the other
hand, from the definition of $\lambda_i^\pm$ $(\ref{UPwind})$ and
the CFL-condition $(\ref{CFL})$ we deduce that
\begin{eqnarray*}
&{}&|\frac{\Delta t}{\Delta
x}\Lambda^+(U_{j}^n-V_{j}^n)|+|(I-\frac{\Delta t}{\Delta
x}\Lambda^++\frac{\Delta t}{\Delta
x}\Lambda^-)(U_{j}^n-V_{j}^n)|+|(-\frac{\Delta
t}{\Delta x}\Lambda^-)(U_{j}^n-V_{j}^n)|\\
&=&\sum_{i=1}^r\frac{\Delta t}{\Delta
x}\lambda_i^+|u_{ij}^n-v_{ij}^n|+\sum_{i=1}^r(1-\frac{\Delta
t}{\Delta x}\lambda_i^++\frac{\Delta t}{\Delta
x}\lambda_i^-)|u_{ij}^n-v_{ij}^n|- \sum_{i=1}^r\frac{\Delta
t}{\Delta x}\lambda_i^-|u_{ij}^n-v_{ij}^n|\\
&=&|U_j^n-V_j^n|.
\end{eqnarray*}
Hence we see that
\begin{eqnarray*}
\sum_{j=-\infty}^{+\infty}|U_j^{n+1} -
V_j^{n+1}|\leq\sum_{j=-\infty}^{+\infty}|U_j^n-V_j^n|\leq\cdots\leq\sum_{j=-\infty}^{+\infty}|U_j^{0}-V_j^{0}|.
\end{eqnarray*}
This completes the proof.
\end{proof}

By taking $V_j^n=U_{j-1}^n$ in Lemma $\ref{lem3}$, we get the
following corollary on BV-estimates of the difference solutions.
\begin{cor}\label{cor2}
Let $U_j^n$ be a solution to the difference scheme $(\ref{dc1})$
with initial data $U_j^0$. Then the $BV$-estimate
\begin{eqnarray*}
 \sum_{j=-\infty}^{+\infty}|U_{j}^{n}-U_{j-1}^{n}|\leq \sum_{j=-\infty}^{+\infty}|U_{j}^{0}-U_{j-1}^{0}|
\end{eqnarray*}
holds for all $n\ge 0$.
\end{cor}

Having Corollary $\ref{cor2}$, we show the uniform boundedness of
the difference solutions, which is not covered in Lemma \ref{lem2}.
\begin{lem}\label{boundedness}
\begin{eqnarray*}
\sup_j|U_j^n|\leq\sum_{j=-\infty}^{+\infty}|U_{j}^{0}-U_{j-1}^{0}|
\end{eqnarray*}
holds for all $n\ge 0$.
\end{lem}

\begin{proof}
Thanks to the initial data $(\ref{inital})$, we have $U_{j_0}^n=0$
if $j_0$ is large enough. Thus, for any $j$ it follows that
\begin{eqnarray*}
|U_j^n|=|U_j^n - U_{j_0}^n|\leq \sum_{k=j + 1}^{j_0}|U_k^n -
U_{k-1}^n|\leq\sum_{k=-\infty}^{+\infty}|U_k^0 - U_{k-1}^0|.
\end{eqnarray*}
This completes the proof.
\end{proof}

\begin{rem}
This lemma requires that the initial data are of bounded variation.
Otherwise, we could not obtain the uniform boundedness of the
difference solutions.
\end{rem}

\section{Time-Lipschitz Continuity}
\setcounter{equation}{0}

In this section we show the time-Lipschitz continuity of the
difference solutions, which seems not so easy for the present
problem.

To begin with, we set
$G(U)=(f_1,f_2,\cdots,f_{r-1})^T\in\mathbb{R}^{r-1}$ and
$G_j^n=G(U_j^n)$. Note that
\begin{eqnarray}\label{QF}
Q(U)=KG(U), \qquad G(U)= K'Q(U),
\end{eqnarray}
where $K$ is the constant $r\times(r-1)$-matrix
\begin{eqnarray}\label{matrix}
K = \left( \begin{array}{cccccc}
1 &  &   &  &    \\
-1 & 1&   & & \\
{}&  -1& 1& & \\
{}& &  \ldots& \ldots &\ldots\\
 {}& &  & -1& 1\\
 {}& &  & & -1
\end{array} \right).
\end{eqnarray}
and $K'$ is the constant $(r-1)\times r$-matrix such that $K'K=I_{r
-1}$. Moreover, it follows from the mean-value theorem and the
definitions in $(\ref{int1})$ and $(\ref{int2})$ that
$$
G_j^{n+1}-G_j^n=\mathcal {M}(U_j^n, U_j^{n+1})(U_j^{n+1}-U_j^n)
$$
with $\mathcal {M}(U_j^n, U_j^{n+1})$ the following $(r-1)\times r$
matrix
\begin{eqnarray*}
\mathcal {M}(U_j^n, U_j^{n+1})=\left( \begin{array}{cccccc}
A_1 & B_1  &   &  &    \\
    &  A_2& B_2  & & \\
{}&  & A_3&B_3 & \\
{}& &  \ldots& \ldots &\ldots\\
 {}& &  & A_{r-1}&B_{r-1}
\end{array} \right).
\end{eqnarray*}
Set $M = M(U_j^n, U_j^{n+1})=\mathcal {M}(U_j^n,
U_j^{n+1})K\in\mathbb{R}^{(r-1)\times (r-1)}$. By calculation, we
have
\begin{eqnarray*}
M =\left( \begin{array}{cccccc}
A_1-B_1 & B_1  &   &  &    \\
  -A_2  & A_2-B_2 & B_2  & & \\
{}&  -A_3& A_3-B_3&B_3 & \\
{}& &  \ldots& \ldots &\ldots\\
 {}& &  & -A_{r-1}&A_{r-1}-B_{r-1}
\end{array} \right).
\end{eqnarray*}
Thanks to the strict monotonicity assumption (1) on the $f_i$'s, it
is an elementary fact that $I-\frac{\Delta t}{\epsilon}M$ is
invertible.

We assume that there exists a $(U_j^n, U_j^{n+1})$-independent norm
$|\cdot|_\star$ on ${\bf R}^{r-1}$ and a positive constant $\lambda$
such that
\begin{eqnarray}\label{inverse}
|(I-\frac{\Delta t}{\epsilon}M)^{-1}|_\star \leq (1 + \frac{\Delta
t}{\epsilon}\lambda)^{-1}.
\end{eqnarray}
Although such a norm has not been found for general cases with $r
> 3$, this assumption is indeed true for $r\le 3$ and for
the case where $Q(U)$ is linear with respect to $U$. In fact, we
have

\begin{prop}\label{prop1}
For $r\le 3$, there exists a positive constant $\lambda$ such that
\begin{eqnarray*}
|(I-\frac{\Delta t}{\epsilon}M)^{-1}|_\infty\leq (1 + \frac{\Delta
t}{\epsilon}\lambda)^{-1}.
\end{eqnarray*}
\end{prop}

\begin{proof}
When $r=2$, $M$ is the negative number $(A_1 - B_1)$. Thus we have
$$
|(I-\frac{\Delta t}{\epsilon}M)^{-1}|_\infty = (1 + \frac{\Delta
t}{\epsilon}(B_1-A_1))^{-1}\leq (1 + \frac{\Delta
t}{\epsilon}\lambda)^{-1}
$$
with $\lambda = \min_{V, W}\{B_1(V, W) - A_1(V, W)\} > 0$. Here the
boundedness of the different solution established in Lemma
\ref{boundedness} has been used.

For $r=3$, $M$ is the $2\times2$ matrix
\begin{eqnarray*}
M =\left( \begin{array}{cccccc}
A_1-B_1 & B_1   \\
-A_{2}&A_{2}-B_{2}
\end{array} \right).
\end{eqnarray*}
Then we have
\begin{eqnarray*}
I-\frac{\Delta t}{\epsilon}M=\left( \begin{array}{cccccc}
1+\frac{\Delta t}{\epsilon}(B_1-A_1) & -\frac{\Delta t}{\epsilon}B_1  \\
\frac{\Delta t}{\epsilon}A_{2}&1+\frac{\Delta
 t}{\epsilon}(B_{2}-A_{2})
\end{array} \right)
\end{eqnarray*}
and
\begin{eqnarray*}
&{}&(I-\frac{\Delta t}{\epsilon}M)^{-1}\\
&=& \frac{1}{1+ A\frac{\Delta t}{\epsilon}+ B\frac{\Delta
 t^2}{\epsilon^2}}
\left( \begin{array}{cccccc} 1+\frac{\Delta t}{\epsilon}(B_2-A_2) &
\frac{\Delta t}{\epsilon}B_1   \\
-\frac{\Delta t}{\epsilon}A_{2}&1+\frac{\Delta t}{\epsilon}(B_1-A_1)
\end{array} \right)
\end{eqnarray*}
with $A=B_1-A_1+B_2-A_2$ and $B=B_1B_2 - A_1B_2 + A_1A_2$. Recall
that $A_i < 0$ and $B_i > 0$ for $i=1, 2$. Set
$$
a=\max\{B_1 + B_2-A_2, B_1-A_1 -A_2\}.
$$
It is obvious that $0 < a < A$ and $B>0$. Thus we deduce that
\begin{eqnarray*}
|(I-\frac{\Delta t}{\epsilon}M)^{-1}|_\infty = \frac{1 +
a\frac{\Delta t}{\epsilon}}{1+ A\frac{\Delta t}{\epsilon}+
B\frac{\Delta
 t^2}{\epsilon^2}}\leq (1 +
\frac{\Delta t}{\epsilon}\lambda)^{-1}
\end{eqnarray*}
with
$$
\lambda = \min\{A - a, B/a\} >0.
$$
This completes the proof.
\end{proof}

For $r>3$, $M$ is no longer strictly diagonally dominant and the
inequality in Proposition \ref{prop1} does not hold anymore.
However, we have

\begin{prop}\label{prop2}
If $Q(U)$ is linear with respect to $U$, then the above assumption
$(\ref{inverse})$ holds.
\end{prop}

\begin{proof}
Since $Q(U)$ is linear with respect to $U$, $M$ is a constant
matrix. Set
$\alpha_1=1,\alpha_{i+1}=-\alpha_{i}\frac{B_i}{A_{i+1}},$ and
$D=diag(\alpha_1,\alpha_2,\cdots,\alpha_{r-1})$. It is obvious that
the diagonal matrix $D$ is positive definite and $DM$ is symmetric.
Then $D^{\frac{1}{2}}MD^{-\frac{1}{2}}$ is also symmetric. Thus,
there exist an orthogonal matrix $T$ and a diagonal
negative-definite matrix $N$ such that
\begin{eqnarray*}
D^{\frac{1}{2}}MD^{-\frac{1}{2}}=T^{-1}N T,
\end{eqnarray*}
and thereby
\begin{eqnarray*}
(TD^{\frac{1}{2}})(I-\frac{\Delta
t}{\epsilon}M)^{-1}=(I-\frac{\Delta
t}{\epsilon}N)^{-1}(TD^{\frac{1}{2}}).
\end{eqnarray*}
Set $P=TD^{\frac{1}{2}}$. We define a norm on ${\bf R}^{r-1}$ as
$$
|\xi|_p = |P\xi|
$$
for $\xi\in {\bf R}^{r-1}$. Thus the corresponding matrix norm is
\begin{eqnarray*}
|(I-\frac{\Delta t}{\epsilon}M)^{-1}|_\star : = \sup_{|\xi|_p\leq
1}|(I-\frac{\Delta t}{\epsilon}M)^{-1}\xi|_p \\
= \sup_{|\xi|_p\leq 1}|P(I-\frac{\Delta
t}{\epsilon}M)^{-1}\xi|\\
= \sup_{|\xi|_p\leq 1}|(I-\frac{\Delta
t}{\epsilon}N)^{-1}P\xi|\\
\leq \sup_{|\xi|_p
\leq 1}|(I-\frac{\Delta t}{\epsilon}N)^{-1}||P\xi|\\
\leq|(I-\frac{\Delta t}{\epsilon}N)^{-1}| = (1 + \frac{\Delta
t}{\epsilon}\lambda)^{-1},
\end{eqnarray*}
where $\lambda$ is the smallest eigenvalue of $-N$.
\end{proof}

Now we turn to estimate $Q^n=Q(U_j^n)$.

\begin{lem}\label{lem4}
Assume (\ref{inverse}) holds. Then there is a positive constant $C$
such that
\begin{eqnarray*}
||Q^n||_{L^1}\le C((1+\lambda\frac{\Delta
t}{\epsilon})^{-n}||Q^0||_{L^1}+\epsilon).
\end{eqnarray*}
 for all $n$.
\end{lem}
\begin{proof}
From the difference scheme it follows that
\begin{eqnarray*}
G_j^{n+1}-G_j^n& = &\mathcal{M}(U_j^n, U_j^{n+1})(U_j^{n+1}-U_j^n)\\
& = & \mathcal{M}(U_j^n, U_j^{n+1})[\frac{\Delta
t}{\epsilon}Q_j^{n+1} - \frac{\Delta t}{\Delta
x}\Lambda^+(U_j^{n}-U_{j-1}^{n}) - \frac{\Delta t}{\Delta
x}\Lambda^-(U_{j+1}^{n}-U_j^{n})]\\
&=&\frac{\Delta t}{\epsilon}MG_j^{n+1}-\mathcal{M}(U_j^n,
U_j^{n+1})[\Lambda^+(U_j^{n}-U_{j-1}^{n})+\Lambda^-(U_{j+1}^{n}-U_j^{n})]\frac{\Delta
t}{\Delta x}
\end{eqnarray*}
and thereby
\begin{eqnarray*}
(I-\frac{\Delta t}{\epsilon}M)G_j^{n+1}=G_j^n-\mathcal {M}(U_j^n,
U_j^{n+1})[\Lambda^+(U_j^{n}-U_{j-1}^{n})\\+
\Lambda^-(U_{j+1}^{n}-U_j^{n})]\frac{\Delta t}{\Delta x}.\nonumber
\end{eqnarray*}
Thus we deduce from the inequality (\ref{inverse}) that
\begin{eqnarray*}
&&(1+\lambda\frac{\Delta t}{\epsilon})|G_j^{n+1}|_\star \leq
|(I-\frac{\Delta t}{\epsilon}M)G_j^{n+1}|_\star \\
& = & |G_j^n-\mathcal {M}(U_j^n,
U_j^{n+1})[\Lambda^+(U_j^{n}-U_{j-1}^{n})+\Lambda^-(U_{j+1}^{n}-U_j^{n})]\frac{\Delta
t}{\Delta x}|_\star\\
& \leq & |G_j^n|_\star + |\mathcal {M}(U_j^n,
U_j^{n+1})[\Lambda^+(U_j^{n}-U_{j-1}^{n}) +
\Lambda^-(U_{j+1}^{n}-U_j^{n})]|_\star\frac{\Delta t}{\Delta x} .
\end{eqnarray*}
Since $|\cdot|_\star$ is equivalent to the $1$-norm $|\cdot|$ on
${\bf R}^{r-1}$ and since the difference solution is bounded, there
is a positive constant $C$ such that
\begin{eqnarray*}
(1+\lambda\frac{\Delta t}{\epsilon})|G_j^{n+1}|_\star \leq
|G_j^n|_\star + C\frac{\Delta t}{\Delta x}(|U_j^{n}-U_{j-1}^{n}| +
|U_{j+1}^{n}-U_j^{n}|).
\end{eqnarray*}
Consequently, we use Corollary $\ref{cor2}$ to obtain
\begin{eqnarray*}
(1+\lambda\frac{\Delta t}{\epsilon})\sum_j|G_j^{n+1}|_\star\Delta x
&\le& \sum_j|G_j^n|_\star\Delta x+C\Delta
t\sum_j|U_j^{n}-U_{j-1}^{n}|\\
&\le& \sum_j|G_j^n|_\star\Delta x + C\Delta t.
\end{eqnarray*}

From the last inequality, it is easy to verify that
\begin{eqnarray*}
\sum_j|G_j^{n}|_\star\Delta x\le (1+\lambda\frac{\Delta
t}{\epsilon})^{-n}\sum_j|G_j^{0}|_\star\Delta x+C\epsilon\le
C((1+\lambda\frac{\Delta t}{\epsilon})^{-n}\sum_j|G_j^{0}|\Delta
x+\epsilon).
\end{eqnarray*}
Thus we get
\begin{eqnarray*}
||G^n||_{L^1} \leq C\sum_j|G_j^{n}|_\star\Delta x \le
C((1+\lambda\frac{\Delta t}{\epsilon})^{-n}\sum_j|G_j^{0}|\Delta
x+\epsilon).
\end{eqnarray*}
Finally, from $(\ref{QF})$ it follows that
\begin{eqnarray*}
||Q^n||_{L^1}\le C||G^n||_{L^1} \qquad \mbox{and} \qquad
||G^0||_{L^1}\le C||Q^0||_{L^1}.
\end{eqnarray*}
Hence
\begin{eqnarray*}
||Q^n||_{L^1}\le C((1+\lambda\frac{\Delta
t}{\epsilon})^{-n}||Q^0||_{L^1}+\epsilon)
\end{eqnarray*}
and the proof is complete.
\end{proof}

Now we can easily show the time-Lipschitz continuity of $U^n$.
\begin{lem}\label{lem5}
Assume (\ref{inverse}) holds. Then
\begin{eqnarray*}
||U_j^{n+1}-U_j^n||_{L^1}\le C(\frac{\Delta
t}{\epsilon}(1+\lambda\frac{\Delta
t}{\epsilon})^{-(n+1)}||Q^0||_{L^1}+\Delta t)
\end{eqnarray*}
for $n=0,1,2,\cdots.$
\end{lem}

\begin{proof}
From the difference scheme $(\ref{dc1})$ we have
\begin{eqnarray*}
U_j^{n+1}-U_j^n=\frac{\Delta t}{\epsilon}Q_j^{n+1}-\frac{\Delta
t}{\Delta x}\Lambda^+(U_j^{n}-U_{j-1}^{n})-\frac{\Delta t}{\Delta
x}\Lambda^-(U_{j+1}^{n}-U_j^{n})
\end{eqnarray*}
and therefore
\begin{eqnarray*}
|U_j^{n+1}-U_j^n|\le\frac{\Delta
t}{\epsilon}|Q_j^{n+1}|+C\frac{\Delta t}{\Delta
x}(|U_j^{n}-U_{j-1}^{n}|+|U_{j+1}^{n}-U_j^{n}|).
\end{eqnarray*}
Thus, from Corollary $\ref{cor2}$ and Lemma $\ref{lem4}$ we deduce
that
\begin{eqnarray*}
||U_j^{n+1}-U_j^n||_{L^1}\le\frac{\Delta
t}{\epsilon}||Q_j^{n+1}||_{L^1}+C\Delta t\leq C(\frac{\Delta
t}{\epsilon}(1+\lambda\frac{\Delta
t}{\epsilon})^{-(n+1)}||Q^0||_{L^1}+\Delta t).
\end{eqnarray*}
This completes the proof.
\end{proof}

\section{Entropies}
\setcounter{equation}{0}

In this section, we discuss entropy functions for the
reaction-hyperbolic system $(\ref{pde1})$ and its equilibrium system
$(\ref{pde2})$. For this purpose, we make the following additional assumption\\

\begin{itemize}

\item[(3).] For any fixed $v$, there exists $w$ such that $vf_i(v, w) \geq
0$ for $i=1, 2, \cdots, r-1$.\\
\end{itemize}

\noindent This assumption, together with those made in the previous
section, ensures that there exists a unique and globally-defined
function $h_i$ of one variable such that
\begin{equation}\label{31}
f_i(u_i, u_{i+1})= 0 \qquad \mbox{iff} \qquad u_{i+1} = h_i(u_i)
\end{equation}
for $i=1, 2, \cdots, r-1$. Obviously, $h_i$ is strictly increasing
and $h_i(0)=0$. By the implicit function theorem, $h_i$ is
continuously differentiable. Thus the equilibrium system
$(\ref{pde2})$ can be written as
\begin{equation}\label{32}
\begin{split}
\partial_t(u_1+u_2+\cdots+u_r) +
\partial_x(\lambda_1u_1+\lambda_2u_2+\cdots+\lambda_ru_r)=0,\\
u_2=h_1(u_1), \ u_3 = h_2(u_2), \ \cdots, \ u_r = h_{r-1}(u_{r-1}).
\end{split}
\end{equation}

Set
\begin{eqnarray*}
v = u_1 + h_1(u_1) + h_2\circ h_1(u_1) + \cdots + h_{r-1}\circ
h_{r-2}\circ \cdots \circ h_1(u_1).
\end{eqnarray*}
Since the right-hand side is strictly increasing with respect to
$u_1$, $u_1$ can be expressed as a function of $v$, say $u_1=
u_1(v)$. By the inverse function theorem, $u_1(v)$ is continuously
differentiable. Set
$$
u_{i+1}(v) = h_i(u_i(v))
$$
for $i=1, 2, \cdots, r-1$. Consequently, each $u_i(v)$ is strictly
increasing, continuously differentiable and $u_i(0)=0$. Set
\begin{eqnarray*}
v=u_1(v)+u_2(v)+\cdots+u_r(v) \quad \mbox{and} \quad h(v) =
\lambda_1u_1(v) + \lambda_2u_2(v) + \cdots + \lambda_r u_r(v).
\end{eqnarray*}
Then the equilibrium system $(\ref{32})$ can be rewritten as
\begin{eqnarray}\label{pde2*}
\partial_tv + \partial_xh(v)=0.
\end{eqnarray}
Recall that any convex function is a convex entropy function for
scalar conservation laws like (\ref{pde2*}) (see, e.g., \cite{Da}).

Next, we turn to discuss the entropy functions for system
$(\ref{pde2*})$ and the reaction-hyperbolic system $(\ref{pde1})$.

\begin{lem}\label{entropy1}
Given a strictly convex smooth function $\tilde{\eta}(v)$, there is
a dissipative entropy function $\eta(U)$, in the sense of \cite{Y3},
for the reaction-hyperbolic system $(\ref{pde1})$ such that
\begin{eqnarray*}
\tilde{\eta}(v)=\eta(u_1(v), u_2(v), \cdots,
u_r(v))\equiv\eta(U(v)).
\end{eqnarray*}
\end{lem}

\begin{proof}
We inductively define
\begin{eqnarray*}
\eta_r(u) &=& \int^u_0\tilde{\eta}'(u_r^{-1}(w))dw, \\
\eta_{i-1}(u) &=& \int^u_0\eta_i'(h_{i-1}(w))dw
\end{eqnarray*}
for $i= r, r-1, \cdots, 2$. Then we have
\begin{eqnarray}\label{35}
\eta_{i-1}'(u_{i-1}) = \eta_i'(h_{i-1}(u_{i-1}))  \qquad \mbox{and}
\qquad \eta_r'(u) = \tilde{\eta}'(u_r^{-1}(u)) .
\end{eqnarray}

Since $\tilde{\eta}(v)$ is strictly convex and $u_r^{-1}$ and $h_i$
are all strictly increasing, it is clear that $\eta_1', \eta_2',
\cdots, \eta_r'$ are strictly increasing. Therefore,
\begin{eqnarray}\label{36}
\eta(U): =\sum_{i=1}^r\eta_i(u_i) + \tilde{\eta}(0)
\end{eqnarray}
is a strictly convex function of $U$. Recall the matrix $K$ defined
in (\ref{matrix}). We see that
$$
(\eta'_1 - \eta'_2, \eta'_2 - \eta'_3, \cdots, \eta'_{r-1} -
\eta'_r)^T = K^T\eta_U(U)
$$
and, moreover, from (\ref{QF}) that
\begin{equation}\label{37}
\begin{split}
Q(U) & = KG(U) \\
&=K\mbox{diag}(\frac{f_1(u_1, u_2)}{\eta'_1 - \eta'_2},
\frac{f_2(u_2, u_3)}{\eta'_2 - \eta'_3}, \cdots,
\frac{f_{r-1}(u_{r-1}, u_r)}{\eta'_{r-1} - \eta'_r})K^T\eta_U(U)\\
&\equiv S(U)\eta_U(U).
\end{split}
\end{equation}
Thanks to the relations in (\ref{35}), we deduce from (\ref{31}),
the convexity of $\eta_i$ and the monotonicity assumption (1) that
\begin{eqnarray*}
\frac{f_i(u_i, u_{i+1})}{\eta'_i - \eta'_{i+1}} &=& \frac{f_i(u_i,
u_{i+1})}{\eta'_{i+1}(h_i(u_i)) -
\eta'_{i+1}(u_{i+1})}\\
& = & -\frac{\int_0^1 f_{iu_{i+1}}(u_i, h_i(u_i) + \sigma(u_{i+1} -
h_i(u_i)))d\sigma}{\int_0^1 \eta''_{i+1}(h_i(u_i) + \sigma(u_{i+1} -
h_i(u_i)))d\sigma} < 0.
\end{eqnarray*}
Thus, $S(U)$ is a symmetric and non-positive definite matrix. Its
null space is obviously that of $K^T$, which is independent of $U$.
Consequently, $\eta(U)$ is a dissipative entropy function, in the
sense of \cite{Y3}, for the reaction-hyperbolic system
$(\ref{pde1})$.

Furthermore, we deduce from the relations in (\ref{35}) that
\begin{eqnarray*}
\eta_1'(u_1(v))&=&\eta_2'(h_1(u_1(v)))\\
&=&\eta_3'(h_2\circ h_1(u_1))\\
&=&\cdots\\
&=&\eta_r'(h_{r-1}\circ h_{r-2}\circ \cdots \circ h_1(u_1))\\
&=&\eta_r'(u_r(v))\\
&=&\tilde{\eta}'(v).
\end{eqnarray*}
Similarly, we have
\begin{eqnarray*}
\eta_2'(u_2(v))=\cdots=\eta_r'(u_r(v))=\tilde{\eta}'(v).\\
\end{eqnarray*}
Therefore, we deduce from $u_1(v)+u_2(v)+\cdots+u_r(v)\equiv v$ that
\begin{eqnarray*}
(\eta(U(v)))'&=&(\sum_{i=1}^r\eta_i(u_i(v)))'\\
&=&\eta_1'(u_1(v))u_1'(v)+\eta_2'(u_2(v))u_2'(v)+\cdots+\eta_r'(u_r(v))u_r'(v)\\
&=&\tilde{\eta}'(v).
\end{eqnarray*}
In view of $\eta_i(0)=0=u_i(0)$, we see from (\ref{36}) that
$\eta(U(0))=\tilde\eta(0)$ and hence
\begin{eqnarray*}
\tilde{\eta}(v)=\eta(U(v)).
\end{eqnarray*}
This completes the proof.
\end{proof}

We conclude this section with a discrete entropy inequality for the
difference solutions.
\begin{lem}\label{entropy}
Let $U_j^n$ be a solution to the difference scheme $(\ref{dc1})$.
Then, for any smooth convex function
$\eta(U)=\sum_{i=1}^r\eta_i(u_i)$, there exists a Lipschitz
continuous function $\Psi$ of two variables such that for all $j\in
\mathbb{Z}$ and $n\geq 0$, the following cell entropy inequalities
hold:
\begin{eqnarray*}
\eta(U_j^{n+1})&\leq& \eta(U_j^n)-\frac{\Delta t}{\Delta
x}(\Psi(U_{j}^n,U_{j+1}^n)-\Psi(U_{j-1}^n,U_{j}^n))+\frac{\Delta
t}{\epsilon}\eta_U(U_j^{n+1})Q(U_j^{n+1}).
\end{eqnarray*}
Moreover, the Lipschitz continuous function satisfies the
consistency relation
\begin{eqnarray*}
\Psi(U,U)=\sum_{i=1}^r\lambda_i\eta_i(u_i).
\end{eqnarray*}
\end{lem}

\begin{proof}
For any smooth convex function $\eta_i$ and any two real numbers
$a,b\in \mathbb{R}$, it is standard that
\begin{eqnarray*}
\eta_i(b)-\eta_i(a)\leq \eta_i'(b)(b-a).
\end{eqnarray*}
Thus, for the given convex function
$\eta(U)=\sum_{i=1}^r\eta_i(u_i)$ and for any $V,W\in \mathbb{R}^r$
we have
\begin{eqnarray*}
\eta(V)-\eta(W)\leq \eta_U(V)(V-W).
\end{eqnarray*}
Thus, it follows from the original difference scheme that
\begin{eqnarray*}
\eta(U_{j}^{n+1})\leq\eta(U_j^n-\frac{\Delta t}{\Delta
x}\Lambda^+(U_j^{n}-U_{j-1}^{n})-\frac{\Delta t}{\Delta
x}\Lambda^-(U_{j+1}^{n}-U_j^{n}))+\eta_U(U_{j}^{n+1})\frac{\Delta
t}{\epsilon}Q_j^{n+1},
\end{eqnarray*}
since
\begin{eqnarray*}
U_j^{n+1}=\frac{\Delta t}{\epsilon}Q_j^{n+1}+U_j^n-\frac{\Delta
t}{\Delta x}\Lambda^+(U_j^{n}-U_{j-1}^{n})-\frac{\Delta t}{\Delta
x}\Lambda^-(U_{j+1}^{n}-U_j^{n}).
\end{eqnarray*}
On the other hand, we deduce from $(\ref{UPwind}), (\ref{CFL})$ and
the convexity of $\eta_i$ that
\begin{eqnarray*}
&{}&\eta(U_j^n-\frac{\Delta t}{\Delta
x}\Lambda^+(U_j^{n}-U_{j-1}^{n})-\frac{\Delta t}{\Delta
x}\Lambda^-(U_{j+1}^{n}-U_j^{n}))\\
&=&\eta(\frac{\Delta t}{\Delta x}\Lambda^+U_{j-1}^n+(I-\frac{\Delta
t}{\Delta x}\Lambda^++\frac{\Delta t}{\Delta
x}\Lambda^-)U_j^n-\frac{\Delta t}{\Delta
x}\Lambda^-U_{j+1}^{n})\\
&=&\sum_{i=1}^r\eta_i(\lambda_i^+\frac{\Delta t}{\Delta
x}u_{i,j-1}^n+(1-\lambda_i^+\frac{\Delta t}{\Delta
x}+\lambda_i^-\frac{\Delta t}{\Delta x})u_{i,j}^n
-\lambda_i^-\frac{\Delta t}{\Delta
x}u_{i,j+1}^n)\\
&\leq&\sum_{i=1}^r[\lambda_i^+\frac{\Delta t}{\Delta
x}\eta_i(u_{i,j-1}^n)+(1-\lambda_i^+\frac{\Delta t}{\Delta
x}+\lambda_i^-\frac{\Delta t}{\Delta
x})\eta_i(u_{i,j}^n)-\lambda_i^-\frac{\Delta t}{\Delta
x}\eta_i(u_{i,j+1}^n)]\\
&=&\sum_{i=1}^r[\eta_i(u_{i,j}^n)-\frac{\Delta t}{\Delta
x}(|\lambda_i|\eta_i(u_{i,j}^n)-\lambda_i^+\eta_i(u_{i,j-1}^n)+\lambda_i^-\eta_i(u_{i,j+1}^n))].
\end{eqnarray*}
Define
\begin{eqnarray*}
\Psi(U,V)=\sum_{i=1}^r[\frac{\lambda_i}{2}(\eta_i(u_i)+\eta_i(v_i))
+\frac{|\lambda_i|}{2}(\eta_i(u_i)-\eta_i(v_i))].
\end{eqnarray*}
This $\Psi$ is obviously Lipschitz continuous and satisfies the
consistency relation. Moreover, the above inequalities lead directly
to
\begin{eqnarray*}
\eta(U_j^{n+1})&\leq& \eta(U_j^n)-\frac{\Delta t}{\Delta
x}(\Psi(U_{j}^n,U_{j+1}^n)-\Psi(U_{j-1}^n,U_{j}^n))+\frac{\Delta
t}{\epsilon}\eta_U(U_j^{n+1})Q(U_j^{n+1}).
\end{eqnarray*}
This completes the proof.
\end{proof}

\section{Main Results}
\setcounter{equation}{0}

In this section we prove the main results of this paper. To this
end, we define
$$
U^\Delta(x,t)=(u_1^{\Delta}(x,t), u_2^{\Delta}(x,t),\cdots,
u_r^{\Delta}(x,t))^T:=(u_{1,j}^{n}, u_{2,j}^{n},\cdots,
u_{r,j}^{n})^T
$$
for $(x,t)\in [j\Delta x,(j+1)\Delta x)\times[n\Delta t,(n+1)\Delta
t)$. With this definition, it simply follows from Corollary
$\ref{cor2}$, Lemma \ref{boundedness} and Lemma \ref{lem5} that
\begin{lem}\label{41}
The piecewise constant function $U^{\Delta}(x,t)$ satisfies the
following estimates
\begin{eqnarray}
&{}&|U^{\Delta}(x,t)|\leq\sum_{j=-\infty}^{+\infty}|U_{j}^{0}-U_{j-1}^{0}| \quad
\textrm{ for all } (x,t),\label{1}\\
&{}&TV(U^{\Delta}(\cdot,t))\leq\sum_{j=-\infty}^{+\infty}|U_{j}^{0}-U_{j-1}^{0}|,\label{estimate1}\\
&{}&||U^\Delta(\cdot,t)-U^\Delta(\cdot,t_1)||_{L^1}\leq
C(\frac{1}{\epsilon}(1+\lambda\frac{\Delta
t}{\epsilon})^{{-\frac{\min\{t,t_1\}}{\Delta t}}}||Q^0||_{L^1}+1)
(|t-t_1|+\Delta t)\label{2}
\end{eqnarray}
for all $t, t_1>0$.
\end{lem}

\begin{proof}
We only need to show the last inequality. Let $k,k_1$ be two
integers such that $t\in [k\Delta t,(k+1)\Delta t),t_1\in [k_1\Delta
t,(k_1+1)\Delta t)$. Without loss of generality, we assume $k_1\leq
k.$ Then we deduce from the definition of $U^\Delta$ and Lemma
$\ref{lem5}$ that
\begin{eqnarray*}
||U^\Delta(\cdot,t)-U^\Delta(\cdot,t_1)||_{L^1}&=&\sum_{j}|U_j^k-U_j^{k_1}|\Delta
x\\&\leq&\sum_{n=k_1}^{k-1}\sum_{j}|U_j^{n+1}-U_j^{n}|\Delta
x\\&\leq&\sum_{n=k_1}^{k-1}C(\frac{1}{\epsilon}(1+\lambda\frac{\Delta
t}{\epsilon})^{-(n+1)}||Q^0||_{L^1}+1)\Delta t\\
&\leq& C(\frac{1}{\epsilon}(1+\lambda\frac{\Delta
t}{\epsilon})^{-(k_1+1)}||Q^0||_{L^1}+1)(k-k_1)\Delta t\\
&\le&C(\frac{1}{\epsilon}(1+\lambda\frac{\Delta
t}{\epsilon})^{-\frac{\min\{t,t_1\}}{\Delta
t}}||Q^0||_{L^1}+1)(|t-t_1| + \Delta t).
\end{eqnarray*}
This completes the proof.
\end{proof}

Having the estimates in Lemma \ref{41} and the discrete entropy
inequality in Lemma $\ref{entropy}$, we follow the standard argument
in \cite{Da} to obtain ($\epsilon$ is fixed)

\begin{thm}\label{th1}
Suppose the initial data
$U_0(x)=(u_{10}(x),u_{20}(x),\cdots,u_{r0}(x))$ have bounded
variations, the grid sizes $\Delta t$ and $\Delta x$ satisfy the
CFL-condition $(\ref{CFL})$, the $f_i's$ satisfy the assumptions
(1)--(3), and the inequality (\ref{inverse}) holds. Then, as the
grid sizes $\Delta t,\Delta x$ tend to zero, there is a subsequence
of the function family $U^{\Delta}(x,t)=(u_1^{\Delta}(x,t),
u_2^{\Delta}(x,t),\cdots, u_r^{\Delta}(x,t))$ converging in
$(L_{loc}^1(\mathbb{R}^1 \times \mathbb{R}^+) )^r$ to an entropy
solution $U^{\epsilon}(x,t)=(u_1^{\epsilon}, u_2^{\epsilon},\cdots,
u_r^{\epsilon})$ of the reaction-hyperbolic system $(\ref{pde1})$
with initial data $U_0(x)$. Furthermore, the solution fulfills the
following estimates
\begin{eqnarray}\label{estimate2}
&{}&|U^{\epsilon}(x,t)|\leq\sum_{j=-\infty}^{+\infty}|U_{j}^{0}-U_{j-1}^{0}|
\quad \textrm{ for almost all } (x,t),\label{3}\\
&{}&TV(U^{\epsilon}(\cdot,t))\leq\sum_{j=-\infty}^{+\infty}|U_{j}^{0}-U_{j-1}^{0}|,\\
&{}&||U^\epsilon(\cdot,t)-U^\epsilon(\cdot,t_1)||_{L^1}\leq
C(\frac{1}{\epsilon}\exp(-\frac{\lambda\min\{t,t_1\}}{\epsilon})||Q^0||_{L^1}+1)|t-t_1|.
\label{4}
\end{eqnarray}
for all $t,t_1>0$.
\end{thm}


In the framework of BV-solutions, the zero-relaxation limit can be
very easily discussed. In fact, the standard argument in \cite{Da}
proves that the embedding of
$L^\infty(\mathbb{R}_+\times\mathbb{R})\cap L^\infty(\mathbb{R}_+,
BV(\mathbb{R}))\cap Lip(\mathbb{R}_+, L^1(\mathbb{R}))$ into
$L^1_{loc}(\mathbb{R}_+\times\mathbb{R})$ is compact. On the other
hand, the estimates in $(\ref{3})-(\ref{4})$ show that
$\{U^\epsilon\}_{\epsilon>0}$ lies in a bounded subset of
$L^\infty(\mathbb{R}_+\times\mathbb{R})\cap L^\infty(\mathbb{R}_+,
BV(\mathbb{R}))\cap Lip(\mathbb{R}_+, L^1(\mathbb{R}))$ by assuming
\begin{equation}\label{47}
\|Q^0\|_{L^1}=0 .
\end{equation}
Namely, the initial data are assumed to be in equilibrium. Thus, we
have
\begin{thm} Under the conditions of
Theorem $\ref{th1}$ and the equilibrium assumption (\ref{47}), there
exist a bounded measurable function
$U^*(x,t)=(u_1^*(x,t),u_2^*(x,t),\cdots,u_r^*(x,t))$ and a
subsequence(denoted in the same way) of set
$\{U^\epsilon(x,t)=(u_1^{\epsilon}, u_2^{\epsilon},\cdots,
u_r^{\epsilon})\}$ such that as $\epsilon\rightarrow0$,
\begin{eqnarray*}
U^\epsilon(x,t) \rightarrow U^*(x,t) &in& (L_{loc}^1(\mathbb{R}^1
\times \mathbb{R}^+) )^r.
\end{eqnarray*}
Moreover, the function $U^*(x,t)$ is a weak entropy solution to the
Cauchy problem $(\ref{pde2})$ with initial data $U_0(x)$ and
satisfies
\begin{eqnarray*}
&{}&|U^{*}(x,t)|\leq\sum_{j=-\infty}^{+\infty}|U_{j}^{0}-U_{j-1}^{0}| \quad \textrm{ for almost all } (x,t),\label{5}\\
&{}&TV(U^{*}(\cdot,t))\leq\sum_{j=-\infty}^{+\infty}|U_{j}^{0}-U_{j-1}^{0}|,\\
&{}&||U^*(\cdot,t)-U^*(\cdot,t_1)||_{L^1}\leq C|t-t_1|.\label{6}
\end{eqnarray*}
\end{thm}

\begin{rem}
The fact that $U^*(x,t)$ satisfies the entropy conditions for the
equilibrium system $(\ref{pde2})$ follows from Lemmas \ref{31} and
\ref{32}. Indeed, because $S(U)$ in (\ref{37}) is symmetric and
non-positive, the term $\eta_U(U^{n+1})Q(U^{n+1})$ in Lemma \ref{32}
is
$$
\eta_U(U^{n+1})Q(U^{n+1})=\eta_U(U^{n+1})S(U^{n+1})\eta_U(U^{n+1})\leq
0
$$
if $\eta$ is chosen to be that constructed in Lemma \ref{31}.
\end{rem}

\begin{rem}
Without the equilibrium assumption in (\ref{47}), $U^\epsilon$ can
only converge to $U^*$ for $t>0$ but not up to $t=0$, because
$$
\exp(-\lambda\min\{t, t_1\}/\epsilon)\leq
\frac{\epsilon}{e\lambda\min\{t, t_1\}}
$$
in $(\ref{4})$. Indeed, without the equilibrium assumption, initial
boundary-layers occur. See also \cite{Y1}.
\end{rem}



\end{document}